\documentclass[reqno]{amsart}
\usepackage{amsmath,amssymb, amsthm}
\usepackage{graphicx}
\usepackage{amsmath,amscd}
\usepackage{mathtools}
\usepackage[margin=1in]{geometry}
\usepackage[mathscr]{euscript}

\usepackage{fancyhdr}

\begin{document}

\newtheorem{problem}{Problem}
\newtheorem{theorem}{Theorem}[section]
\newtheorem{definition}{Definition}[section]
\newtheorem{lemma}[theorem]{Lemma}

\title
{A Converse Theorem for Metaplectic Eisenstein Series on $SL_2(\mathbb{A})$}

\author{ Vladislav Petkov}
\maketitle

\begin{abstract}The purpose of this work is to produce a converse theorem for adelic Eisenstein series on the double metaplectic cover of the group $SL_2(\mathbb{A})$. We show that the double Dirichlet series, which satisfy the natural functional equations required in our theorem, are in fact Mellin transforms at infinity of metaplectic Eisenstein series.\end{abstract}

\begin{section}{Introduction}

Following Hecke's original result, number theorists have produced various converse theorems for classical groups. Perhaps the most celebrated such result is that of Cogdell and Piatetski-Shapiro \cite{Cogdell-PS}, where they give a criterion to check the automorphicity of irreducible cuspidal representations of $GL_n(\mathbb{A})$. To this day, however, few have investigated the analogous question for the $n-$fold metaplectic covers of the classical groups.

The method developed in \cite{Cogdell-PS} can be adopted for the cases when strong multiplicity one applies. Gelbart and Piatetski-Shapiro named these representations "exceptional" or "distinguished" and in particular proved in \cite{GelbartPS2} that, in the case of the double metaplectic cover of the group $SL_2(\mathbb{A})$, these are precisely the so called "elementary theta series". Nevertheless, following this approach one faces two important obstructions. First, the existence of cuspidal distinguished representations for higher metaplectic covers or higher rank $n$ is still unclear, although they are conjectured to correspond to classical cuspidal representations of lower rank (see \cite{Suzuki}). And more importantly a uniqueness of the Whittaker model is highly unusual for general metaplectic representations. An alternative method will be to translate the problem via the Shimura or higher theta correspondences and use the classical non-metaplectic results.

 On the other hand, it is viable to consider the continuous spectrum of the metaplectic group. The first to do so were Godlfeld and Diamantis in \cite{GoldfeldDiamantis1}. Later in \cite{GoldfeldDiamantis2} they showed that their converse theorem can be used to associate double zeta functions coming from prehomogeneous vector spaces to the Fourier-Whittaker coefficients of the classical\footnote{Throughout this paper the word "classical" is used both for "non-adelic" and "non-metaplectic". We trust that in each context there will not be a confusion.} metaplectic Eisenstein series and, as a result, prove extra functional equations for those double zeta functions. These hidden equations will come from the  relation between Whittaker coefficients of metaplectic Eisenstein series and Weyl group multiple Dirichlet series that was established by  Brubaker, Bump and Friedberg, \cite{EisensteinWeylMult}.

 The purpose of this paper is to generalize their result and give a criterion to relate double Dirichlet series to arbitrary adelic metaplectic Eisentein series on $SL_2(\mathbb{A})$. For simplicity we work over the field $\mathbb{Q}$ and henceforth $\mathbb{A}$ will mean the adele ring $\mathbb{A}_\mathbb{Q}$. However, we believe that the proof will work over an arbitrary number field of class number one.
 We would like to mention that our method can be adjusted to give an alternative converse theorem for the standard non-metaplectic Eisenstein series.

To ease the reader we summarize here the layout of this work. In the next Section \ref{Notation} we review some basic notation and properties of the metaplectic group and Eisenstein series. Sections \ref{SecTwists} and \ref{SectionLfunctions} are dedicated to the twists of the metaplectic Eisenstein series and the associated $L-$functions respectively. In these two sections we introduce the main tools we need for our proof. In Section \ref{SectionMainResults} we define a "nice" family of double Dirichlet series, motivated by the previous sections, and prove the Converse theorem that is the center piece of this article. Finally in Section \ref{SectionConclusion} we make some remarks for of the possible applications and generalizations of this result.

To put our theorem in context recall that the double metaplectic cover of $G(\mathbb{A})=SL_2(\mathbb{A})$ is a central extension

$$1\rightarrow \mu_2\rightarrow \overline{G}(\mathbb{A})\rightarrow G(\mathbb{A})\rightarrow 1$$

We define the group law of $\overline{G}(\mathbb{A})$ as

$$(g_1,\epsilon_1)(g_2,\epsilon_2):=(g_1g_2,\epsilon_1\epsilon_2\sigma(g_1,g_2)).$$

Above $g_i\in G(\mathbb{A})$, $\epsilon_i\in \mu_2$ and $\sigma(\cdot,\cdot)$ is a particular 2-cocycle, explicitly defined in \cite{GelbartBook}. We return to this definition more carefully in the next Section. Since there is a natural inclusion of $G\hookrightarrow \overline{G}$ we usually use $g$ as a shorthand for $(g,1)$.

Recall that the \emph{genuine} automorphic forms and representations of $\overline{G}$ are those not coming from trivial lifts of forms and representations of $G$.

Let $\chi_\mathbb{A}\,:\, \mathbb{A}/\mathbb{Q}\rightarrow \mathbb{C}$ be a Hecke character of conductor $N'$. We will see that when we consider genuine automorphic forms on $\overline{G}(\mathbb{A})$ this conductor must be divisible by $4$. Let $I(\bar{g},w)$ be  a smooth function induced from the action of $\chi_\mathbb{A}$ on the Levi subgroup, such that
\begin{equation}\label{Idef0}
         I\left(\left[\left(
                        \begin{array}{cc}
                          a & b \\
                          0 & a^{-1} \\
                        \end{array}
                      \right)
         ,\epsilon\right]\bar{g},w\right)=\epsilon\cdot\chi_\mathbb{A}(a)|a|^{w}I(\bar{g},w),
\end{equation}

for all $a\in \mathbb{A}^\times$, $b\in \mathbb{A}$, $\epsilon\in \mu_2$, $\bar{g}\in \overline{G}(\mathbb{A})$ and $w\in \mathbb{C}$.

If $P\in G$ is the subgroup of upper triangular matrices, the standard Eisenstein series associated to the function $I(\bar{g},w)$ is defined as

$$E(\bar{g},w,I):= \sum_{\gamma\in P(\mathbb{Q})\backslash G(\mathbb{Q})} I(\gamma \bar{g}, w)$$

Now we are prepared to state the main result.

\begin{theorem}[The Converse Theorem]\label{ConverseTheorem} Let $\chi_\mathbb{A}:\mathbb{A}^\times/\mathbb{Q}^\times\rightarrow \mathbb{C}$ be the Hecke character above. Let $D$ range over the reduced set of residues modulo $N'=4N$ and let $\chi$ range over the set of Dirichlet characters modulo $D$.
Let $\tau(\chi)$ be the usual Gauss sum, given by

$$\tau_n(\chi):=\sum_{\begin{array}{c}  m\,(\textrm{mod}\, D) \\
                                                                        (m,D)=1
                                                                      \end{array} } \chi(m)e^{2\pi i mn/D}.$$

Let $\mathcal{F}$ be a family of double Dirichlet series

$$L^\pm_j(s,w,\chi):=\sum_{\pm n>0}\sum_{m>0}\frac{a^j_{n,m}\tau_n(\chi)}{m^w|n|^s},$$

where the coefficients $a^j_{n,m}$ are complex numbers and the index $j=1,\ldots, m_N$, for a constant $m_N$, determined by $N$. Assume that this family is "nice" and satisfies the conditions in Definition \ref{Nice Family}. Then there exists a function $I(\bar{g},w)$, as defined in (\ref{Idef0}), so that the double Dirichlet series in the above family are Mellin transforms at infinity of the metaplectic Eisenstein series $E(\bar{g},w,I)$.

\end{theorem}

\end{section}

\begin{section}{Metaplectic Eisenstein series and basic notation}\label{Notation}

In this section we recall some basic terminology and notation for metaplectic Eisenstein series.

Let $p$ be a place of $\mathbb{Q}$ (finite or infinite) and let $\overline{G}_{p}$ be the two fold central extension on $G_{p}=SL(2, \mathbb{Q}_p)$ by $\mu_2=\pm1$.

$$1\rightarrow \mu_2\rightarrow \overline{G}_{p}\rightarrow SL_2(\mathbb{Q}_p)\rightarrow 1.$$

Such non-trivial extension exists and is unique, since the place $p$ is not complex. We define the group law of $\overline{G}$ as

$$(g_1,\epsilon_1)(g_2,\epsilon_2):=(g_1g_2,\epsilon_1\epsilon_2\sigma_p(g_1,g_2)).$$

Above $g_i\in G_{p}$, $\epsilon_i\in \pm 1$ and $\sigma_p(\cdot,\cdot)$ is the local component of the 2-cocycle $\sigma(\cdot ,\cdot)$(ref. to \cite{GelbartBook}). To ease notation we will write $g=(g,1)\in \overline{G}_{p}$, utilizing the natural inclusion of $G_{p} \hookrightarrow \overline{G}_{p}$.

Let $N_p$, $M_p$ and $K_p$ denote respectively the unipotent, Levi and maximal compact subgroups of $G$. In particular

$$N_p=\left\{ n(b)= \left(
                \begin{array}{cc}
                  1 & b \\
                  0 & 1 \\
                \end{array}
              \right)\,:\, b\in\mathbb{Q}_p
\right\}
,\,\,M_p=\left\{ m(a)=\left(
                        \begin{array}{cc}
                          a & 0 \\
                          0 & a^{-1} \\
                        \end{array}\right)\,:\, a\in\mathbb{Q}_p^\times
 \right\}$$

Recall $K_p=SO(2,\mathbb{R})$, if $p=\infty$, and $K_p=SL(2,\mathbb{Z}_p)$, if $p<\infty$.

Let $P_p=N_pM_p$ be the standard Borel subgroup. By the Iwasawa decomposition we have $G_p=N_pM_pK_p$.

Let $\overline{N}_p$, $\overline{M}_p$, $\overline{K}_p$, $\overline{P}_p=\overline{N}_p\overline{M}_p$ be the corresponding metaplectic lifts inside $\overline{G}$.

Define the \emph{global} metaplectic group as
$$\overline{G}_\mathbb{A}=\left.\left(\prod  \nolimits ^\prime \overline{G}_{p}\right)\right/ Z_0 ,$$

Where $Z_0= \left\{ \epsilon=(\epsilon_p)\in \otimes_p \mu_{2} | \prod_p \epsilon_p =1 \right\}$.

Let $\overline{N}_\mathbb{A}$, $\overline{M}_\mathbb{A}$ and $\overline{K}_\mathbb{A}$ be the global equivalents of $\overline{N}_p$, $\overline{M}_p$, $\overline{K}_p$.

Fix once and for all an additive character $\psi$ of the unipotent subgroup $N_\mathbb{A}$, thinking of it as a character of $\mathbb{A}/\mathbb{Q}$, and give $N_\mathbb{A}$ the normalized self-dual measure with respect to $\psi$, i.e. $vol(N_\mathbb{Q}\backslash N_\mathbb{A})=1$.

Let $\chi_\mathbb{A}:\, \mathbb{A}^\times/\mathbb{Q}^\times\rightarrow \mathbb{C}$ be an idele class character, to which we associate an induced genuine representation $\mathcal{I}(w,\chi_\mathbb{A})$, defined as the space of smooth functions $I(\bar{g},w)$ on $\overline{G}(\mathbb{A})$, such that

\begin{equation}\label{Idef}
         I((p,\epsilon)\bar{g},w)=\epsilon\cdot\chi_\mathbb{A}(a)|a|^{w}I(\bar{g},w),
\end{equation}
where $p=n(b)m(a)$. Note that we use a normalization different from the one in \cite{GelbartBook}(p. 65, (3.16)).

We will consider only $I(\cdot,w)$ that are factorizable, i.e. $I(\cdot,w)=\otimes I_p(\cdot,w)$, where $I_p(\cdot,w)$ come from the corresponding local representations. The standard Eisenstein series associated to $I(\cdot,s)$ is defined as

\begin{equation}
E(\bar{g},w,I):= \sum_{\gamma\in P(\mathbb{Q})\backslash G(\mathbb{Q})} I(\gamma\bar{g},w).
\end{equation}

This series is convergent for $Re(w)>2$, has a meromorphic continuation to $w\in\mathbb{C}$ and satisfies the functional equation
\begin{equation}
E(\bar{g},1-w,M(w)I)=E(\bar{g},w,I),
\end{equation}

where $M(w): \mathcal{I}(w,\chi_\mathbb{A})\rightarrow \mathcal{I}(1-w,\chi_\mathbb{A}^{-1})$ is the intertwining operator, defined in the half-plane of absolute convergence by

$$M(w)I(\bar{g},w)=\int_\mathbb{A} I(wn(b)\bar{g},w)db,$$
for $w=\left(
         \begin{array}{cc}
           0 & -1 \\
           1 & 0 \\
         \end{array}
       \right)
$. The Eisenstein series has a Fourier expansion with respect to the character $\psi$

\begin{equation}
E(\bar{g},w,I)=\sum_{m\in\mathbb{Q}}E_m(\bar{g},w,I)
\end{equation}
where
\begin{equation}
E_m(\bar{g},w,I)=\int_{\mathbb{Q}\backslash \mathbb{A}} E(n(b)\bar{g},w,I)\psi(-mb)db.
\end{equation}

When $I(\cdot,w)=\otimes I_p(\cdot,w)$ is factorizable the Fourier coefficients have product expansions

$$E_m(\bar{g},w,I)=\prod_{p\leq \infty}W_{m,p}(\bar{g}_p,w,I_p)$$
where the local Whittaker function is defined as
$$W_{m,p}(\bar{g}_p,w,I_p)= \int_{\mathbb{Q}_p}I_p(wn(b)\bar{g}_p,w)\psi_p(-mb)db.$$

In order to prove a converse theorem for the Eisenstein series $E(\bar{g},w,I)$ we  need to somewhat modify it. Let $\bar{g}=(g,\epsilon)$ and let $g=(g_\infty,g_f)$, with $g_\infty$ and $g_f$ being the infinite and finite components of $g$ respectively. Write $$g_\infty =\left(
                                  \begin{array}{cc}
                                    1 & x \\
                                    0 & 1 \\
                                  \end{array}
                                \right)\left(\begin{array}{cc}
                                         y^{1/2} & 0 \\
                                         0 & y^{-1/2}
                                       \end{array}\right)\left(
                                                           \begin{array}{cc}
                                                             \cos(\theta) & \sin(\theta) \\
                                                             -\sin(\theta) & \cos(\theta) \\
                                                           \end{array}
                                                         \right)
                                       $$

Setting $z=x+iy$ and $\bar{g}_z=((g_\infty,1),1)$, we write $$E(z,w,I)=E(\bar{g}_z,w,I).$$

From (\ref{Idef}) and strong approximation $E(\bar{g},w,I)$ is uniquely determined by its values on the set $\{\bar{g}_z: z\in \mathfrak{H}\}$. As a consequence $E(\bar{g},w,I)$ and its Fourier coefficients are also determined by the series $E(z,w,I)$, which we will consider instead.

The series $E(z,w,I)$ looks like a modular form of the upper half-plane, although it is not immediately clear of which weight or level. These are determined by the properties of $I$ as follows.

Define a character of $\overline{K}_\infty$ as $$\nu_l(\bar{k})=\nu_l((k_\theta,\epsilon)):=(\epsilon)^{2l} e^{il\theta},$$ for $l\in \frac{1}{2}\mathbb{Z}$, such that $$I_\infty(\bar{g}\bar{k},w)=\nu_l(\bar{k})I(\bar{g},w).$$

Note that $\nu_l$ depends only on $l\, (\textrm{mod} \,\,2)$.

Let $S$ be the set of finite primes where $I$ is not \emph{spherical}, i.e. it is not constant on $K_p$ for $p\in S$. These are exactly the places where $\chi_\mathbb{A}$ is ramified. Let $p^{t_p}\mathbb{Z}_p$ be the conductor of $\chi_\mathbb{A}$ at each $p\in S$ and let $$N':=\prod_{p\in S}p^{t_p}.$$

Note that since we chose the representation $\mathcal{I}(w,\chi)$ to be \emph{genuine},  $2$ is always bad prime, so $N'=4N$ and $l\in \frac{1}{2}+\mathbb{Z}$. Then the Eisenstein series $E(z,w,I)$ is a modular form of weight $l$ for the congruence group $\Gamma_0(4N)$.

Recall that the space of non-holomorphic Eisenstein series for a modular group $\Gamma$ is generated by the Eisenstein series at the cusps of $\Gamma$. In order to define those, we remind the reader of some classical notation used in \cite{Roelke} and \cite{GoldfeldDiamantis2}. Let $\nu(\gamma)$ be the weight $l$ multiplier system for $SL(2,\mathbb{R})$, coming from $\chi_\mathbb{A}$ above.  For $\gamma=\left(
                                                                                                                          \begin{array}{cc}
                                                                                                                            a &b \\
                                                                                                                            c &d \\
                                                                                                                          \end{array}
                                                                                                                        \right)
$ let $j(\gamma,z)$ be the automorphic factor, defined as

\begin{equation}
j(\gamma,z)=\frac{(cz+d)^{l/2}}{(c\bar{z}+d)^{l/2}},
\end{equation}

and let the slash operator to be defined for $f:\frak{H}\rightarrow \mathbb{C}$ and $\gamma\in SL(2,\mathbb{R})$ as follows

\begin{equation}\label{automorphiceqn}
(f|\gamma)(z)=j(\gamma,z)^{-1}f(\gamma z).
\end{equation}

This operator satisfies the following relation:

\begin{equation}\label{rfactor}
f|M|N=r(M,N)f|(MN),\,\, (M,N\in SL(2,\mathbb{R})),\end{equation}

where

\begin{equation}
r(M,N):=\frac{(c_MNz+d_M)^{l/2}(c_Nz+d_N)^{l/2}}{(c_{MN}z+d_{MN})^{l/2}},\,\, (\text{ for } M,N\in SL_2(\mathbb{R}))
\end{equation}

In fact, the following lemma holds.
\begin{lemma}\label{Rfactordef}
Let $M=\left(
         \begin{array}{cc}
           \ast & \ast \\
           m_1 & m_2 \\
         \end{array}
       \right)
$, $N=\left(
         \begin{array}{cc}
           a & b \\
           c & d \\
         \end{array}
       \right)
$ and $MN=\left(
         \begin{array}{cc}
           \ast & \ast \\
           m_1' & m_2' \\
         \end{array}
       \right)
$, where

$M, N \in SL(2,\mathbb{R})$ are any matrices. Then $r(M,N)= e^{\frac{\pi i l}{2}w(M,N)}$, with

$$w(M,N)=\left\{\begin{array}{ll}
                  sgn(c)+sgn(m_1)-sgn(m_1')-sgn(m_1c m_1')& m_1c m_1'\neq0, \\
                  (sgn(c)-1 )(1-sgn(m_1))& m_1c\neq 0,\, m_1'=0,\\
                  (sgn(c)+1 )(1-sgn(m_2))& m_1'c\neq 0, \,m_1=0,\\
                  (1-sgn(a))(1+sgn(m_1))& m_1m_1'\neq0, \, c=0, \\
                  (1-sgn(a))(1-sgn(m_2))& m_1=m_1'=c=0.
                \end{array}
   \right. $$

\end{lemma}

The proof of the lemma follows easily from the definition of $r(M,N)$ in \cite{Roelke},(\S 1.) and Theorem 16 in \cite{Maass}.

Next let $\{\frak{a}_i: 1\leq i\leq m_N\}$ be the set of inequivalent singular cusps of $\Gamma_0(4N)$,i.e. those for which $\nu(\gamma_{\frak{a}_i})=1$, where $\gamma_{\frak{a}_i}$ is the generator of the stabilizer group $\Gamma_{\frak{a}_i}\subset \Gamma_0(4N)$ of the cusp $\frak{a}_i$. We do not need to consider the other cusps, since by \cite{Roelke} (10.7) the only Eisenstein series that matter are those at the singular cusps. For convenience we choose the ordering, so that $\frak{a}_1=\infty$ and $\frak{a}_{m_N}=0$. Let $\sigma_{\frak{a}_i}(\infty)=\frak{a}_i$ be the scaling matrix, such that $\sigma_{\frak{a}_i}^{-1} \Gamma_{\frak{a}_i}   \sigma_{\frak{a}_i}=\Gamma_\infty$. In particular, from above $\sigma_{\frak{a}_1}=I$ and $\sigma_{\frak{a}_{m_N}}=W_{4N}$, where $I$ is the identity and $W_{4N}$ is the Fricke involution $\left(
              \begin{array}{cc}
                0 & -1/(2\sqrt{N}) \\
                2\sqrt{N} & 0 \\
              \end{array}
            \right).
$

Finally we are ready to define the Eisenstein series at the cusp $\frak{a}_i$, for $1\leq i \leq m_N$. For $Re(w)>1$ they are given by the following convergent series

\begin{equation}
E_i(z,w):= \sum_{\gamma\in\Gamma_{\frak{a}_i}\backslash \Gamma_0(4N)}\nu(\gamma)^{-1}j(\gamma,z)^{-1}r(\sigma_{\mathfrak{a}_i}^{-1},\gamma)^{-1}Im(\sigma_{\frak{a}_i}\gamma z)^w
\end{equation}

The following lemma summarizes the properties of the metaplectic Eisenstein series that we will need through this paper. For a concise and beautiful proof of all statements we refer the reader to \cite{Roelke}, \S 10.

\begin{lemma}
Let $\frak{a}_j$ be any cusp, $j=\frak{a}_1,\ldots,\frak{a}_{m_N}$. The Eisenstein series satisfies the following properties

\begin{itemize}
\item[(a)] $E_i(\cdot, w)|\gamma=\nu(\gamma)E_i(\cdot,w)$, for any $\gamma\in \Gamma_0(4N)$.
\item[(b)] The Eisenstein series has meromorphic continuation to $w\in\mathbb{C}$.
\item[(c)] Let $E(z,w):=(E_1(z,w),\ldots,E_{m_N}(z,w))^T$ then
\begin{equation}\label{scateringeqn}
E(z,1-w)=\Phi(1-w)E(z,w)
\end{equation}
\item[(d)] $E_i(z,w)$ is an eigenfunction of $\triangle_l$ with eigenvalue $w(1-w)$.

\item[(e)] For any pair of  cusps there is a Fourier expansion of $E_i(z,w)$ at $\frak{a}_j$

\begin{equation}
E_i(z,w)|\sigma_{\frak{a}_j}=\delta_{ij}y^w+p_{ij}(w)y^{1-w} + \sum_{n\neq 0} a^{ij}_n(w) W_{\frac{sgn(n)k}{2},w-\frac{1}{2}}(4\pi|n|y)e^{2\pi i n x}.
\end{equation}
\end{itemize}

\end{lemma}

In the above notation  $\triangle_l$ is the weight $l$ Laplacian:

$$\triangle_l:= -y^2\left(\frac{\partial^2}{\partial x^2} +\frac{\partial^2}{\partial y^2}\right)+ ily\frac{\partial}{\partial x}$$ $\Phi(w)=\left(p_{ij}(w)\right)$ is the scattering matrix(see \cite{Roelke}) and $W_{a,b}$ is the classical Whittaker function (\cite{Whittaker})

\begin{equation}
W_{a,b}(z)=\frac{e^{-\frac{z}{2}} z^a}{\Gamma(\frac{1}{2}-a+b)}\int_0^\infty u^{-a-\frac{1}{2}+b}(1+z^{-1}u)^{a-\frac{1}{2}+b}e^udu.
\end{equation}

\end{section}

\begin{section}{Twists of Eisenstien series} \label{SecTwists}

Let $N$ be as in the previous section. For every $D$, such that $(D,4N)=1$ let $\chi$ be a Dirichlet character modulo $D$. For the purpose of the converse theorem, that is the main result of this work, we
need to introduce the twists of our Eisenstein series by these Dirichlet characters. For any function $f:\, \frak{H}\times \mathbb{C}\rightarrow \mathbb{C}$ let

\begin{equation}\label{twistedfunction}
f(\cdot,\chi):=\sum_{\begin{array}{c}  m\,(\textrm{mod}\, D )\\
                                                                        (m,D)=1
                                                                      \end{array} }\chi(m)f|\left(
                                                                                             \begin{array}{cc}
                                                                                               1 & m/D \\
                                                                                               0 & 1 \\
                                                                                             \end{array}
                                                                                           \right)
                                                                      \end{equation}

Recall that the Gauss sum $\tau_n(\chi)$ is defined as

$$\tau_n(\chi):=\sum_{\begin{array}{c}  m\,(\textrm{mod}\, D) \\
                                                                        (m,D)=1
                                                                      \end{array} } \chi(m)e^{2\pi i mn/D}.$$

When the Fourier expansion of $f(z,w)$ at infinity is

\begin{equation}
f(z,w)=a(w)y^w+b(w)y^{1-w} + \sum_{n\neq 0} a_n(w) W_{\frac{sgn(n)l}{2},w-\frac{1}{2}}(4\pi|n|y)e^{2\pi i n x},
\end{equation}

the corresponding expansion of the twist will be

\begin{equation}\label{TwistedExpansion}
f(z,w,\chi)=\tau_0(\chi)(a(w)y^w+b(w)y^{1-w}) + \sum_{n\neq 0}\tau_n(\chi) a_n(w) W_{\frac{sgn(n)l}{2},w-\frac{1}{2}}(4\pi|n|y)e^{2\pi i n x}.
\end{equation}

We introduce the following convenient function

\begin{equation}\label{Checkfunction}
\check{f}(z,w):=e^{\frac{\pi i l}{2}}f(z,w)|W_{4N}.
\end{equation}

Observe that $\check{\check{f}}=f$ and $\check{f}(iy)=f(i/(4Ny))$. Next consider $f(z,w,\chi)|W_{4ND^2}$:

\begin{equation}
f(z,w,\chi)\left|\left(
              \begin{array}{cc}
                0 & \frac{-1}{2D\sqrt{N}} \\
                2D\sqrt{N} & 0 \\
              \end{array}
            \right)\right. = \sum_{\begin{array}{c}  m \,(\textrm{mod}\, D) \\
                                                                        (m,D)=1
                                                                      \end{array} }\chi(m)f\left| \left(
                                                                                             \begin{array}{cc}
                                                                                               1 & m/D \\
                                                                                               0 & 1 \\
                                                                                             \end{array}
                                                                                           \right)\right. \left| \left(
              \begin{array}{cc}
                0 & \frac{-1}{2D\sqrt{N}} \\
                2D\sqrt{N} & 0 \\
              \end{array}
            \right)\right.
\end{equation}

Using (\ref{rfactor}) we get

\begin{equation}
f(\cdot,\chi)\left|W_{4ND^2}\right. = \sum_{\begin{array}{c}  m\,(\textrm{mod}\, D)  \\
                                                                        (m,D)=1
                                                                      \end{array} }\chi(m)r\left( \left(
                                                                                             \begin{array}{cc}
                                                                                               1 & m/D \\
                                                                                               0 & 1 \\
                                                                                             \end{array}
                                                                                           \right), \left(
              \begin{array}{cc}
                0 & \frac{-1}{2D\sqrt{N}} \\
                2D\sqrt{N} & 0 \\
              \end{array}
            \right)\right)f\left| \left(
                                                                                             \begin{array}{cc}
                                                                                               2m\sqrt{N} & \frac{-1}{2D\sqrt{N}} \\
                                                                                               2D\sqrt{N} & 0 \\
                                                                                             \end{array}
                                                                                           \right)\right.
\end{equation}

Next using twice Lemma \ref{Rfactordef}, the explicit formula for $w(\cdot,\cdot)$ and (\ref{Checkfunction}), we obtain that the above equals

\begin{eqnarray*}
\lefteqn{e^{-\pi i l/2}\sum_{\begin{array}{c}  m \,(\textrm{mod}\, D)  \\
                                                                        (m,D)=1
                                                                      \end{array} }\chi(m)\check{f}\left| \left(
                                                                                             \begin{array}{cc}
                                                                                               0 & \frac{-1}{2\sqrt{N}} \\
                                                                                               2\sqrt{N} & 0 \\
                                                                                             \end{array}
                                                                                           \right) \left(
              \begin{array}{cc}
                                                                                               2m\sqrt{N} & \frac{-1}{2D\sqrt{N}} \\
                                                                                               2D\sqrt{N} & 0 \\
                                                                                             \end{array}
            \right)\right.=}\\
            && e^{-\pi i l/2}\sum_{\begin{array}{c}  m \,(\textrm{mod}\, D)  \\
                                                                        (m,D)=1
                                                                      \end{array} }\chi(m)\check{f}\left| \left(
                                                                                             \begin{array}{cc}
                                                                                               D & 0 \\
                                                                                               -4mN & 1/D \\
                                                                                             \end{array}
                                                                                           \right)\right.
\end{eqnarray*}

Let $r,t$ be integers, such that $Dt-4Nmr=1$. As $m$ runs over the reduced set of residues $(\textrm{mod}\, D) $, so does $r$ and, as a result, we can write

\begin{equation}\label{Equation3}
f(z,w,\chi)|{W_{4ND^2}}=e^{\frac{-\pi i l}{2}}\overline{\chi(-4N)}\sum_{\begin{array}{c}
                                                                        r\,(\textrm{mod}\, D)  \\
                                                                        (r,D)=1
                                                                      \end{array}
 } \overline{\chi(r)} \check{f}\left|\left(
              \begin{array}{cc}
                D & -r \\
                -4mN & t \\
              \end{array}
            \right)\right.
\left|\left(
      \begin{array}{cc}
        1 & r/D \\
        0 & 1 \\
      \end{array}
    \right)
\right. .
\end{equation}

Following the approach in \cite{GoldfeldDiamantis1}, \cite{GoldfeldDiamantis2} we want to define a new Dirichlet character $\check{\chi}$ modulo $D$, in order to obtain an equation of this type:

$$f(z,w,\chi)|W_{4ND^2}=e^{\frac{-\pi i l}{2}}\overline{\chi(-4N)}H(D,N)\check{f}(z,w,\check{\chi}),$$
where $H(D,N)$ will be some number depending on $D$ and $N$, that comes from the multiplier system $\nu$.

In \cite{GoldfeldDiamantis1} they compute $H(D,N)$ exactly for the weight $l=\frac{1}{2}$ and $\nu=\nu_\vartheta$, where $\nu_\vartheta$ is the \emph{theta multiplier system}. For the convenience of the reader we restate the definition of $\nu_\vartheta$ in the form given in \cite{GoldfeldDiamantis2}. Let $\gamma=\left(
                                                                              \begin{array}{cc}
                                                                                a & b \\
                                                                                c & d \\
                                                                              \end{array}
                                                                            \right)
\in \Gamma_0(4N)$ and let $\left(\frac{c}{d}\right)$ be the Kronecker symbol.  The theta multiplier system is defined as

\begin{equation}
\nu_\vartheta:=\left( \frac{c}{d}\right)\left\{\begin{array}{ll}
                                  1 \; &d\equiv 1\, (\textrm{mod}\, 4) \\
                                  i \; &d\equiv 3\,(\textrm{mod}\, 4) .
                                \end{array}\right.
\end{equation}

 For a thorough discussion of $\nu_\vartheta$ we refer the reader either to \cite{Knopp}(Chapter 4) or \cite{Maass} (p. 98). Since the purpose of this paper is to generalize the results in \cite{GoldfeldDiamantis2} to any weight $l$ and any multiplier system $\nu$ we need to find a more general definition of $\check{\chi}$ and $H(D,N)$.

To do this we make the following observations. First, $\nu_l$ depends only on $l \,(\textrm{mod}\, 2) $. Second, in $\nu_1$ and $\nu_2$ prime are multiplier systems of weight $k\in \mathbb{R}$ then $\nu_o=\nu_1/\nu_2$ is a multiplier system of weight $0$, hence is an abelian character of $\Gamma_0(4N)$.  Maass shows in \cite{Maass}, that for any congruence subgroup $\Gamma$ there are precisely $|\Gamma/K_*|$ different multiplier systems of each weight $k\in \mathbb{R}$, where $K_*= \langle[\Gamma,\Gamma], -I\rangle$. He also shows how to compute the $6$ different systems for the simplest case $\Gamma(1)=SL(2,\mathbb{Z})$. In the case of $\Gamma_0(4N)$ such computation will be more involved and, in order to ease our notation we will only consider one multiplier system per weight $l$. This will not affect the methods in our proof, since a change by an abelian character $\nu_0$ will only slightly change the factors in the functional equations of the twisted Eisenstein series.

Thus it remains only to find a multiplier system $\nu_l$ of weight $l= \frac{3}{2}$. To do this note that if $\nu_1$ and $\nu_2$ are of weights $k_1$ and $k_2$ respectively, $\nu_1\nu_2$ is a multiplier of weight $k=k_1+k_2$. Thus, $\nu_l\nu_\vartheta$ will be of weight $2$ and hence will be trivial for our purposes. As a result we can simply chose
$$\nu_l=\left\{\begin{array}{c}
          \nu_\vartheta, \, l\equiv \frac{1}{2} \,(\textrm{mod}\, 2) , \\
          \nu_\vartheta^{-1}, \, l\equiv \frac{3}{2} \,(\textrm{mod}\, 2) .
        \end{array}\right.
$$

From the computations in \cite{GoldfeldDiamantis1}, we see that in either case we can define the same Dirichlet character

\begin{equation}
\check{\chi}:=\left(\frac{r}{D}\right)\overline{\chi(r)},
\end{equation}
where $\left(\frac{r}{D}\right)$ is the Jacobi symbol. Let $\epsilon_D$ be given by

$$\epsilon_D=\left\{
 \begin{array}{ll}
                       1 & D\equiv 1 \mod 4,\\
                       -i & D\equiv 3\mod 4, \, l\equiv \frac{1}{2} \,(\textrm{mod}\, 2) ,\\
                       i & D\equiv 3\mod 4, \, l\equiv \frac{3}{2} \,(\textrm{mod}\, 2) .
                     \end{array}
\right.$$

As a result from (\ref{Equation3}) we get

\begin{equation}\label{Equation4}
f(z,w,\chi)|{W_{4ND^2}}=e^{\frac{-\pi i l}{2}}\overline{\chi(-4N)}\left(\frac{4N}{D}\right)\epsilon_D\sum_{\begin{array}{c}
                                                                        r\mod\, D \\
                                                                        (r,D)=1
                                                                      \end{array}
 } \check{\chi}(r) \check{f}\left|\left(
      \begin{array}{cc}
        1 & r/D \\
        0 & 1 \\
      \end{array}
    \right)
\right. .
\end{equation}

Thus, with (\ref{twistedfunction}) and (\ref{Checkfunction}) the following holds for all $1\leq j \leq m_N$

\begin{equation}\label{Equation5}
E_j(z,w,\chi)|W_{4ND^2}=e^{\frac{-\pi i l}{2}}\overline{\chi(-4N)}\left(\frac{4N}{D}\right)\epsilon_D\check{E}_j(z,w,\check{\chi}).
\end{equation}

Recall that $W_{4ND^2}$ takes $\frak{a}_1=\infty$ to $\frak{a}_{m_N}=0$ and thus the main term of the twisted series is

\begin{equation}
\check{a}_0(z,w,\chi)=\overline{\chi(-4N)}\left(\frac{4N}{D}\right)\epsilon_D \tau_0(\check{\chi})\left(\delta_{jm_N}y^w+p_{jm_N}(w)y^{1-w}\right)
\end{equation}

\begin{section}{The associated $L-$functions}\label{SectionLfunctions}

Let $E_j(z,w)$ be the Eisenstein series at the cusp $\frak{a}_j$ and let $a_n^{j}(w)$ be the $n-$th coefficient of its Fourier expansion at $\frak{a}_1=\infty$. We associate to it the following (untwisted) $L-$functions

\begin{equation}
L^{\pm}_j(s,w):= \sum_{\pm n>0}\frac{a_n^j(w)}{|n|^s}.
\end{equation}

Equivalently, for a Dirichlet character $\chi\, \mod D$, we define

\begin{equation}
L^{\pm}_j(s,w,\chi):= \sum_{\pm n>0}\frac{\tau_n(\chi)a_n^j(w)}{|n|^s}.
\end{equation}

We also define the "completed" $L-$functions\footnote{Note that we use the notation in \cite{GoldfeldDiamantis2}.}:

\begin{equation}
\Lambda_j(s,w,u,\chi):=\int_0^\infty \left(E_j((i+u)y, w,\chi)-\tau_0(\delta_{j1}y^w+p_{j1}(w)y^{1-w})\right)y^s\frac{dy}{y}.
\end{equation}

Let $\check{L}_j$ and $\check{\Lambda}_j$ be the corresponding functions for $\check{E}_j(z,w,\check{\chi})$.

For $u\in \mathbb{R}$, $Re(s)\gg1$ and $Re(w)\gg1$, we have from \cite{FunctionsList} (13.23.4)

\begin{equation}\label{LambdatoL}
\Lambda_j(s,w,u,\chi)=\textbf{c}(s,w,u)\left(L^{+}_j(s,w,\chi),L^{-}_j(s,w,\chi)\right)^T,
\end{equation}

where

\begin{eqnarray*}
\lefteqn{
\textbf{c}(s,w;u)= \frac{\Gamma(w+s)\Gamma(s-w+1)}{(4\pi)^s}}\\&&\cdot\left( \frac{F(s-w+1,s+w,s+1-\frac{l}{2};\frac{1+iu}{2})}{\Gamma(s+1-\frac{l}{2})},
\frac{F(s-w+1,s+w,s+1+\frac{l}{2};\frac{1-iu}{2})}{\Gamma(s+1+\frac{l}{2})}\right).
\end{eqnarray*}

Above $F(a,b,c;d)$ is the Gaussian hypergeometric function (see. \cite{FunctionsList}).

Define also
$$L^\pm_E(s,w,\chi):=(L^\pm_1(s,w,\chi),\ldots,L^\pm_{m_N}(s,w,\chi))^T,$$
$$\Lambda_E(s,w,u,\chi):=(\Lambda_1(s,w,u,\chi),\ldots,\Lambda_{m_N}(s,w,u,\chi))^T.$$

Next we will derive the functional equation of the completed $L-$function.

 If we substitute $(i-u)/(2\sqrt{N}D(u^2+1)y)$ in (\ref{Equation5}),
we obtain

\begin{equation}
E_j\left(\frac{(u+i)y}{2\sqrt{N}D}, w; \chi\right)=\overline{\chi(-4N)}\left(\frac{4N}{D}\right)\epsilon_D\left(\frac{1+iu}{1-iu}\right)^{l/2}\check{E}_j\left(\frac{i-u}{2\sqrt{N}D(u^2+1)y},w,\check{\chi}\right).
\end{equation}
Above we used the elementary equality

$$(u+i)^{l}|u+i|^{-l}=e^{\pi i l/2}(1-iu)^{l/2}(1+iu)^{-l/2}.$$

Set $H=\overline{\chi(-4N)}\left(\frac{4N}{D}\right)\epsilon_D\left(\frac{1+iu}{1-iu}\right)^{l/2}$. Following the standard Riemann computation we get:

\begin{eqnarray}\label{RiemmansTrick}
\nonumber\lefteqn{(2\sqrt{N}D)^{-s}\Lambda_j(s,w,u,\chi)=}\\
&&\nonumber= \int_{\frac{1}{\sqrt{u^2+1}}}^\infty \left[E_j\left(\frac{(u+i)y}{2\sqrt{N}D}, w ,\chi\right)-\tau_0(\chi)\left(\delta_{j1}\left(\frac{y}{2\sqrt{N}D}\right)^w+p_{j1}(w)\left(\frac{y}{2\sqrt{N}D}\right)^{1-w}\right)\right]y^s\frac{dy}{y}\\
&&\nonumber +\int_0^{\frac{1}{\sqrt{u^2+1}}} \left[H\cdot\check{E}_j\left(\frac{i-u}{2\sqrt{N}D(u^2+1)y}, w ,\check{\chi}\right)-\tau_0(\chi)\left(\delta_{j1}\left(\frac{y}{2\sqrt{N}D}\right)^w+p_{j1}(w)\left(\frac{y}{2\sqrt{N}D}\right)^{1-w}\right)\right]y^s\frac{dy}{y}\\
&&\nonumber=\int_{\frac{1}{\sqrt{u^2+1}}}^\infty \left\{
\left[E_j\left(\frac{(u+i)y}{2\sqrt{N}D}, w ,\chi\right)-\tau_0(\chi)\left(\delta_{j1}\left(\frac{y}{2\sqrt{N}D}\right)^w+p_{j1}(w)\left(\frac{y}{2\sqrt{N}D}\right)^{1-w}\right)\right]y^s\right.\\
&&\nonumber+H\cdot\left[\check{E}_j\left(\frac{(u+i)y}{2\sqrt{N}D}, w ,\check{\chi}\right)\right.
\\\nonumber&&\left.\left.-e^{\frac{\pi i l }{2}}\tau_0(\check{\chi})\left(\delta_{jm_N}\left(\frac{y}{2\sqrt{N}D}\right)^w+p_{jm_N}(w)\left(\frac{y}{2\sqrt{N}D}\right)^{1-w}\right)
\right](y(u^2+1))^{-s}\right\} \frac{dy}{y}\\
&&\nonumber+(u^2+1)^{-s}\left[(2\sqrt{N}D)^{-w}(u^2+1)^{\frac{s-w}{2}}\left(H\cdot\tau_0(\check{\chi})e^{\frac{\pi i l}{2}}\frac{\delta_{jm_N}}{s-w}-\tau_0(\chi)\frac{\delta_{j1}}{s+w}\right)\right.\\
&&\left. +(2\sqrt{N}D)^{w-1}(u^2+1)^{\frac{s+w-1}{2}}\left(H\cdot\tau_0(\check{\chi})e^{\frac{\pi i l}{2}}\frac{p_{jm_N}(w)}{w+s-1}-\tau_0(\chi)\frac{p_{j1}(w)}{s-w+1}\right)\right].
\end{eqnarray}

The integrals above converge to an entire function in $s$, since the Whittaker function decays exponentially as $y\rightarrow \infty$. Representing $\Lambda_j(s,w,u,\chi)$ in such form allows us to prove the following lemma.

\begin{lemma}\label{PropertiesLambda}

The completed $L-$functions satisfy the following
\begin{itemize}
\item[1.] The function $\Lambda_j(s,w,u,\chi)$ has a meromorphic continuation to $(s,w)\in\mathbb{C}^2$.
\item[2.] If we modify $\Lambda_j(s,w,u,\chi)$ removing the last summand in the above expression we get an entire function, bounded on vertical strips (EBV).
\item[3.] For all $j$ the function $\Lambda_j(s,w,u,\chi)$ satisfies the following equation
\begin{equation}
\Lambda_j(s,w,u,\chi)=H \cdot (4ND^2)^{-s}(1+u^2)^{-s}\check{\Lambda}_j(-s,w,-u,\check{\chi}).
\end{equation}
\item[4.] For $w$, such that neither $w$ nor $1-w$ is a pole of the scattering matrix $\Phi(w)$, we have
\begin{equation}
\Lambda_E(s,1-w,u,\chi)=\Phi(1-w)\Lambda_E(s,w,u,\chi).
\end{equation}
\end{itemize}
\end{lemma}


The above properties are derived from equations (\ref{RiemmansTrick}) and (\ref{scateringeqn}). Note that using (\ref{LambdatoL}) we can change the last property to
\begin{equation}
L_E^\pm(s,1-w,\chi)=\Phi(1-w)L_E^\pm(s,w,\chi).
\end{equation}
To prove this observation, we can just use the linear independence of the hypergeometric functions in the components of $\textbf{c}$, seen as functions in the variable $u$.
\end{section}

\begin{section}{Main Result}\label{SectionMainResults}
In this section we state and prove the converse theorem that is the main result of this work. Following the steps of \cite{GoldfeldDiamantis2}, we define a \emph{"nice family"} of double Dirichlet series as follows.

\begin{definition}\label{Nice Family}
Fix an integer $N\geq 1$ and let $l$, $m_N$, $\nu(\gamma)$ be as in the previous sections. For each $j=1,\ldots, m_N$, let $\{a^j_{n,m}:n\in \mathbb{Z},\,m\in\mathbb{N}\}$ be a sequence of complex numbers. Assume that $|a^j_{n,m}|= O(|n|^{\alpha},m^{\beta})$ have polynomial growth, as $|n|,m\rightarrow \infty$. Let $D$ range over the reduced set of residues modulo $4N$ and let $\chi$ range over the set of Dirichlet characters modulo $D$. For a pair of complex variables $(s,w)$, such that $Re(w)$ and $Re(s)$ are sufficiently large, define a family of Dirichlet series $\mathcal{F}(N,l,\nu)$ as

$$L^\pm_j(s,w,\chi):=\sum_{\pm n>0}\sum_{m>0}\frac{a^j_{n,m}\tau_n(\chi)}{m^w|n|^s}.$$

Let $\Lambda_j(s,w,u,\chi)$ be defined by the equation (\ref{LambdatoL}). Let $L^\pm(s,w,\chi)=(L^\pm_1(s,w\chi),\ldots,L^\pm_{m_N}(s,w,\chi))^T$.
Recall that $H:=\overline{\chi(-4N)}\left(\frac{4N}{D}\right)\epsilon_D\left(\frac{1+iu}{1-iu}\right)^{l/2}$.

We will say that $\mathcal{F}(N,l,\nu)$ is "nice", if there exist a dual family $\check{\mathcal{F}}(N,l,\nu)$

$$\check{L}^\pm_j(s,w,\chi):=\sum_{\pm n>0}\sum_{m>0}\frac{\check{a}^j_{n,m}\tau_n(\chi)}{m^w|n|^s},$$

satisfying the following properties.

\begin{itemize}
\item[(A)] For all $L^\pm_j(s,w,\chi)\in \mathcal{F}(N,l,\nu)$ the completed series $\Lambda_j(s,w,u,\chi)$ have meromorphic continuation to $\mathbb{C}^2$.
    \item[(B)]There exist meromorphic functions $a_j(w), b_j(w),\check{a}_j(w), \check{b}_j(w)$, with no poles in some right half-plane, so that the following are EBV for $Re(w)\gg 1$.

        \begin{eqnarray}
        \nonumber\lefteqn{(2\sqrt{N}D)^s\Lambda_j(s,w,u,\chi)}\\
        &&\nonumber-(u^2+1)^{-s}\left[(2\sqrt{N}D)^{-w}(u^2+1)^{\frac{s-w}{2}}\left(H\cdot\tau_0(\check{\chi})\frac{\check{b}_j(w)}{s-w}-\tau_0(\chi)\frac{b_j(w)}{s+w}\right)\right.\\
&&\left.+(2\sqrt{N}D)^{w-1}(u^2+1)^{\frac{s+w-1}{2}}\left(H\cdot\tau_0(\check{\chi})\frac{\check{a}_j(w)}{s+w-1}-\tau_0(\chi)\frac{a_j(w)}{s-w+1}\right)\right]
        \end{eqnarray}

        \item[(C)]For all $j$ we have
        \begin{equation}
        \Lambda_j(s,w,u,\chi)=H \cdot (4ND^2)^{-s}(1+u^2)^{-s}\check{\Lambda}_j(-s,w,-u,\check{\chi}).
        \end{equation}
        \item[(D)] Let $\Phi(w)$ be the scattering matrix as before. Then the $L-$functions satisfy the functional equations

            \begin{equation}
            L^\pm(s,1-w,\chi)=\Phi(1-w)L^\pm(s,w,\chi).
            \end{equation}

            \item[(E)] Define

\begin{equation}
a_n^j(w)=\sum_{m=1}^\infty \frac{a^j_{n,m}}{m^w}, \,\,\,\, \check{a}_n^j(w)=\sum_{m=1}^\infty \frac{\check{a}^j_{n,m}}{m^w}.
\end{equation}
Assume that, for each fixed $j$ and $w$, with $Re(w)$ large enough, $|a^j_n(w)|, |\check{a}^j_n(w)|=O(|n|^\alpha)$ have polynomial growth, as $|n|\rightarrow \infty$.
\end{itemize}

\end{definition}

We are now ready to prove that, under these assumptions, a "nice family" of double Dirichlet series can be associated to the Mellin transforms at infinity of metaplectic Eisenstein series.

\begin{proof}

Define, for $j=1,\dots m_N$ and $z \in \frak{H}$

\begin{equation}
f_j(z,w)=a_j(w)y^w+b_j(w)y^{1-w}+\sum_{n\neq0}a_n^j(w)W_{\frac{sgn(n)l}{2},w-\frac{1}{2}}(4\pi|n|y)e^{2\pi i n x},
\end{equation}
where $a_j(w), b_j(w)$ are the meromorphic functions from property $(B)$.

As in every converse theorem like those in \cite{GoldfeldDiamantis2},\cite{Cogdell-PS}, we begin by proving the invariance of $f_j(z,w)$ under the modular group $\Gamma_0(4N)$.

As in Section \ref{SecTwists}, for every character $\chi\,(\textrm{mod}\, D)$, $Re(w)\gg 1$, $u\in \mathbb{R}$ and $y>0$, let

\begin{equation}
F_j((i+u)y,w,\chi):=\sum_{n\neq0}a^j_n(w)\tau_n(\chi)W_{\frac{sgn(n)l}{2},w-\frac{1}{2}}(4\pi|n|y)e^{2\pi n uy},
\end{equation}

\begin{equation}
\check{F}_j((i+u)y,w,\chi):=\sum_{n\neq0}\check{a}^j_n(w)\tau(\chi)W_{\frac{sgn(n)l}{2},w-\frac{1}{2}}(4\pi|n|y)e^{2\pi n uy}.
\end{equation}

From the exponential decay of the Whittaker functions and $|a^j_n(w)|, |\check{a}^j_n(w)|=O(|n|^\alpha)$ the Mellin integrals converge absolutely uniformly and thus

\begin{equation}\label{Mellintransform}
\int_0^\infty F_j((i+u)y,w,\chi)y^s\frac{dy}{y}=\Lambda_j(s,w,u,\chi),\,\,\int_0^\infty \check{F}_j((i+u)y,w,\chi)y^s\frac{dy}{y}=\check{\Lambda}_j(s,w,u,\chi).
\end{equation}

When we fix $Re(w)$, the hypergeometric factors in the tensor $\textbf{c}(s,w,u)$ decay exponentially for $|s|\rightarrow \infty$, provided $Re(s)$ is also sufficiently big and $|u|<\delta$ in some small neighbourhood of zero. Thus, through Mellin inversion we get

\begin{equation}
F_j((i+u)y,w,\chi)=\frac{1}{2\pi i}\int_{Re(s)=\sigma_0}\Lambda_j(s,w,u,\chi)y^{-s}ds,
\end{equation}

\begin{equation}
\check{F}_j((i+u)y,w,\chi)=\frac{1}{2\pi i}\int_{Re(s)=\sigma_0}\check{\Lambda}_j(s,w,u,\chi)y^{-s}ds.
\end{equation}

As usual we need to pick $\sigma _0$ large enough for the integrals to converge. From the above remark of the hypergeometric factors we see that by Phragm\`{e}n-Lindel\"{o}f principle we may shift the line of integration to $\sigma_1=-\sigma_0$. Thus we obtain

\begin{eqnarray}\label{Equation6}
F_j((i+u)y,w,\chi)=\frac{1}{2\pi i}\int_{Re(s)=\sigma_1}\Lambda_j(s,w,u,\chi)y^{-s}d s+\sum_{s_0}Res_{s=s_0}\Lambda_j(s,w,u,\chi)y^{-s}.
\end{eqnarray}

An easy computation gives

\begin{eqnarray}\label{Equation7}
\nonumber\lefteqn{\sum_{s_0}Res_{s=s_0}\Lambda_j(s,w,u,\chi)y^{-s}=H\cdot\tau_0(\check{\chi})\left(\check{b}_j(w)(4ND^2(1+u^2)y)^{-w}\right.}\\
&&\left.+\check{a}_j(w)(4ND^2(1+u^2)y)^{w-1}\right)-\tau_0(\chi)\left(b_j(w)y^w+a_j(w)y^{1-w}\right).
\end{eqnarray}

Using assumption (C), we see that the integral in (\ref{Equation6}) equals

\begin{eqnarray}\label{Equation8}
\nonumber\lefteqn{\int_{Re(s)=\sigma_1}H \cdot (4ND^2)^{-s}(1+u^2)^{-s}\check{\Lambda}_j(-s,w,-u,\check{\chi})y^{-s}d s}\\&&=
\int_{Re(s)=\sigma_0}H \cdot (4ND^2)^{s}(1+u^2)^{s}\check{\Lambda}_j(s,w,-u,\check{\chi})y^{s}d s.
\end{eqnarray}

Define the twists of $f_j(z,w,\chi)$ by the twisted Fourier expansion in (\ref{TwistedExpansion}). Then we have

\begin{equation}
f_j((i+u)y,w,\chi)=\tau_0(\chi)\left(b_j(w)y^w+a_j(w)y^{1-w}\right) + F_j((i+u)y,w,u,\chi),
\end{equation}

\begin{equation}
\check{f}_j((i+u)y,w,\check{\chi})=\tau_0(\check{\chi})\left(b_j(w)y^w+a_j(w)y^{1-w} \right)+\check{F}_j((i+u)y,w,u,\check{\chi}),
\end{equation}

Thus, by (\ref{Equation6}), (\ref{Equation7}), (\ref{Equation8}) we get

\begin{equation}
f_j((i+u)y,w,\chi)=H\cdot \check{f}_j \left(\frac{i-u}{4ND^2(1+u^2)y},w,\check{\chi}\right).
\end{equation}
Since this holds for every $z=uy+iy$ we get as in (\ref{Equation5})

\begin{equation}
f_j(z,w,\chi)|W_{4ND^2}=e^{\frac{-\pi i l}{2}}\overline{\chi(-4N)}\left(\frac{4N}{D}\right)\epsilon_D\check{f}_j(z,w,\check{\chi}).
\end{equation}

Finally, by (\ref{Equation3}) and (\ref{Equation4}) we have

\begin{equation}
\sum_{\begin{array}{c}
                                                                        r\mod\, D \\
                                                                        (r,D)=1
                                                                      \end{array}
 } \overline{\chi(r)} \check{f}_j\left|\left(
              \begin{array}{cc}
                D & -r \\
                -4mN & t \\
              \end{array}
            \right)\right.
\left|\left(
      \begin{array}{cc}
        1 & r/D \\
        0 & 1 \\
      \end{array}
    \right)
\right. = \sum_{\begin{array}{c}
                                                                        r\mod\, D \\
                                                                        (r,D)=1
                                                                      \end{array}
 } \overline{\chi(r)}\left(\frac{4Nr}{D}\right)\epsilon_D \check{f}_j
\left|\left(
      \begin{array}{cc}
        1 & r/D \\
        0 & 1 \\
      \end{array}
    \right)
\right. .
\end{equation}

If we use character summation, as in \cite{GoldfeldDiamantis2}, we get

\begin{equation}
\check{f}_j\left|\left(
              \begin{array}{cc}
                D & -r \\
                -4mN & t \\
              \end{array}
            \right)\right.=\left(\frac{4Nr}{D}\right)\epsilon_D \check{f}_j.
\end{equation}
Equivalently

\begin{equation}
f_j\left|\left(
              \begin{array}{cc}
                t & m \\
                4Nr & D \\
              \end{array}
            \right)\right.=\left(\frac{4Nr}{D}\right)\epsilon_D f_j.
\end{equation}

Since these matrices generate $\Gamma_0(4N)$ (see for example \cite{GoldfeldDiamantis1}) we have proven that $f_j(z,w)$ is $\Gamma_0(4N)$ invariant. Proving the moderate growth of $f_j(z,w)$ at every cusp can be done, by using the growth of the coefficients $a^j_{n,m}$ and the elementary bound for the Whittaker function

$$W_{\frac{sgn(n)l}{2},w-\frac{1}{2}}(z)=e^{-\frac{z}{2}}z^{\frac{sgn(n)l}{2}}(1+O(|z|^{-1})).$$

For more details on this part of the proof look in \cite{GoldfeldDiamantis1}.

Therefore, by the spectral decomposition $f_j(z,w)$ must be a linear combination of Maass forms, Eisenstein series and residues of Eisenstein series of the same weight $l$ and multiplier system $\nu$ \footnote{Note that $\nu$ comes into the Converse theorem via the inclusion of the symbol $\epsilon_D$ in the functional equations for the Dirichlet series.}. We still have to remove the possibility for a non-trivial contribution from the discrete spectrum. However, this is not a problem since, if $M(z,w)$ is a component in $f_j(z,w)$, it must be a modular form with eigenvalue $w(1-w)$ for all $Re(w)\gg 1$. However, the discrete spectrum of $\triangle_l$ is in the region $[0,\infty)$, (see \cite{Roelke}), and as a result we get a contradiction when $Re(w)$ is large enough.

Therefore, $f_j(z,w)$ is a linear combination of the Eisenstein series defined in Section \ref{Notation} and by (\ref{Mellintransform}) and (\ref{LambdatoL}) we see that the "nice" family of double Dirichlet series realizes as the sums of Mellin transforms of metaplectic Eisenstein series.

Define $f(z,w)=(f_1(z,w),\ldots f_{m_N}(z,w))^T$. Then there exists a $m_N\times m_N$ matrix $A(w)$ of functions in $w$, which allows us to rewrite the above statement in more compact way
\begin{equation}\label{Matrix A}
f(z,w)=A(w)E(z,w),
\end{equation}
where $E(z,w)=(E_1(z,w)\ldots,E_{m_N}(z,w))$ is as in Section \ref{Notation}. Further, if $A(w)$ is meromorphic and $w$, $1-w$ are not poles of $\Phi(w)$, we have also the following equation

\begin{equation}\label{A-Phi equation}
\Phi(1-w)A(w)\Phi(w)=A(1-w).
\end{equation}

If $A(w)$ is meromorphic the function $f(z,w)$ has a meromorphic continuation to all $w\in \mathbb{C}$.

Finally, the proof of equation (\ref{A-Phi equation}) follows from the functional equation (\ref{scateringeqn}) of $E(z,w)$.

Let $f_{\text{adelic}}(\bar{g},w)$ be the tensor of adelic lifts of the automorphic functions $f_j(z,w)$. Then by the spectral decomposition $f_{\text{adelic}}(\bar{g},w)$ will be a tensor of adelic metaplectic Eisenstein series $$E(\bar{g},w, I)=\left(E(\bar{g},w, I_1),\ldots,E(\bar{g},w, I_{m_N})\right).$$

Note that the different functions $I_j$ may not be  pure tensor products of local factors. However, all appearing pure tensors $\otimes I_{j,p}(\cdot,w)$ will correspond to the same multiplier system $\nu_l$ and a character $\chi_\mathbb{A}$ as defined in Section \ref{Notation}.

\end{proof}

Note that if we compute the actual entries of the matrix $A(w)$ we will be able to completely determine each function $I_j(\bar{g},w)$. This could be very useful and give us a much greater information on the actual arithmetic properties of  the coefficients of our double Dirichlet series.

\end{section}

\end{section}

\begin{section}{Conclusion}\label{SectionConclusion}

The method of using the metaplectic Eisenstein series at infinity is not new. In fact the meromorphic continuation and functional equation of the adelic analog are derived through such computations, \cite{GelbartSally}.
Considering the series on the upper half-plane is essential as it "reveals" the interesting arithmetic information hidden in the adelic Whittaker coefficients of $E(\bar{g},w,I)$. Recently, it was used in \cite{KudlaYang} to compute the values of the series and their first derivatives at critical points.

In \cite{GoldfeldDiamantis2} they relate a certain pair of Shintani zeta functions coming from a prehomogeneous vector space to the metaplectic Eisenstein series for $\Gamma_0(4N)$. They have remarked this could mean they have to satisfy certain hidden functional equations. Furthermore, using their "scalar" converse theorem they give an example where the components of $A(w)$ can be computed.

 In such case we can work our way backwards in Section \ref{Notation} and determine the exact adelic series $E(\bar{g},w,I)$ from the components of the series $A(w)E(z,w)$, by characterizing $I(\bar{g},w)$ completely. As a result, we believe that by  computing explicitly the Whittaker coefficients like in \cite{KudlaYang}, we can not only prove the extra functional equations coming from Weyl group multiple Dirichlet series, but gain a much deeper arithmetic information for the coefficients of our double Dirichlet series.

Although we work over the field of rational numbers it is very viable to adopt our method to an arbitrary number field of class number one. The only impeding factor is the complexity of the calculations, involving the multiplier systems at each infinite place. Alternatively one can try to extend the Converse theorem to higher rank groups or higher order metaplectic coverings. Such theorems will have even more promissing applications to connect multiple zeta functions coming from prehomogeneous  vector spaces and Weyl group multiple Dirichlet series, which appear as Whittaker coefficients of metaplectic Eisenstein series.
\end{section}

\end{document}